\documentclass{amsart}

\theoremstyle{plain}
\newtheorem{thm}{Theorem}[section]

\newtheorem{example}{Example}[section]
\theoremstyle{definition}

\newtheorem{defi}[thm]{Definition}

\begin{document}

\author{Ewa Graczy\'{n}ska}
\address{Opole University of Technology, Institute of Mathematics
\newline
ul. Luboszycka 3, 45-036 Opole, Poland}
\email{e.graczynska@po.opole.pl \hspace{8mm} www.egracz.po.opole.pl}

\title{Dependence spaces}
\maketitle



\begin{abstract}
The Steinitz exchange lemma is a basic theorem in linear algebra
used, for example, to show that any two bases for a
finite-dimensional vector space have the same number of elements.
The result is named after the German mathematician Ernst Steinitz.

We present here another proof of the result of N.J.S. Hughes
\cite{1} on Steinitz' exchange theorem for infinite bases. In our
proof we  assume Kuratowski-Zorn Maximum Principle instead of well
ordering. We present some examples of dependence spaces of general
nature with theirs possible applications of the result in other as
linear or universal algebra domains of mathematical sciences. The
lecture was presented on 77th Workshop on General Algebra, 24th
Conference for Young Algebraists in Potsdam (Germany) on 21st March
2009.
\end{abstract}

\section{Introduction on Steinitz}

We recall some facts of Article by: J.J. O'Connor and E.F. Robertson
from an internet biography of Steinitz from School of Mathematics
and Statistics University of St Andrews, Scotland. The URL of this
page
is:\\
http://www-history.mcs.st-andrews.ac.uk/Biographies/Steinitz.html .

Born: 13 June 1871 in Laurahütte, Silesia, Germany (now Huta Laura,
Poland)

Died: 29 Sept 1928 in Kiel, Germany.

Ernst Steinitz entered the University of Breslau in 1890. He went to
Berlin to study mathematics there in 1891 and, after spending two
years in Berlin, he returned to Breslau in 1893. In the following
year Steinitz submitted his doctoral thesis to Breslau and, the
following year, he was appointed Privatdozent at the Technische
Hochschule Berlin - Charlottenburg.

The offer of a professorship at the Technical College of Breslau saw
him return to Breslau in 1910. Ten years later he moved to Kiel
where he was appointed to the chair of mathematics at the University
of Kiel.

Steinitz was a friend of Toeplitz. The direction of his mathematics
was also much influenced by Heinrich Weber and by Hensel's results
on p-adic numbers in 1899. In \cite{5} interesting results by
Steinitz are discussed. These results were given by Steinitz in
1900, when he was a Privatdozent at the Technische Hochschule Berlin
- Charlottenburg, at the annual meeting of the Deutsche
Mathematiker-Vereinigung in Aachen. In his talk Steinitz introduced
an algebra over the ring of integers whose base elements are
isomorphism classes of finite abelian groups. Today this is known as
the Hall algebra. Steinitz made a number of conjectures which were
later proved by Hall.

Steinitz is most famous for work which he published in 1910. He gave
the first abstract definition of a field in Algebraische Theorie der
Körper in that year. Prime fields, separable elements and the degree
of transcendence of an extension field are all introduced in this
1910 paper. He proved that every field has an algebraically closed
extension field, perhaps his most important single theorem.

The now standard construction of the rationals as equivalence
classes of pairs of integers under the equivalence relation: (a, b)
is equivalent to (c, d) if and only if ad = bc was also given by
Steinitz in 1910.

Steinitz also worked on polyhedra and his manuscript on the topic
was edited by Rademacher in 1934 after his death.

For more information on the role of Steinitz papers consult the book
chapter 400 Jahre Moderne Algebra, of \cite{5}.

\section{Notation}

We use the notation of \cite{1}: $a,b,c,...,x,y,z,...$ (with or
without suffices) to denote the elements of ${\bf S}$ and
$A,B,C,...,X,Y,Z,...$ for subsets of ${\bf S}$, $\mathbb{X}$,
$\mathbb{Y}$,... denote a family of subsets of ${\bf S}$, $n$ is
always a positive integer.

$A+B$ denotes the union of sets $A$ and $B$, $A-B$ denotes the
difference of $A$ and $B$, i.e.is the set of those elements of $A$
which are not in $B$.

\section{Dependent and independent sets}

The following definition is due to N.J.S. Hughes, invented in 1962
in \cite{1}:
\begin{defi}
A set ${\bf S}$ is called a {\it dependence space} if there is
defined a set $\Delta$, whose members are finite subsets of ${\bf
S}$, each containing at least 2 elements, and if the Transitivity
Axiom is satisfied.
\end{defi}

\begin{defi}
A set $A$ is called {\it directly dependent}if $A \in \Delta$.
\end{defi}

\begin{defi}
An element $x$ is called {\it dependent on A} and is denoted by $x
\sim \Sigma A$ if either $x \in A$ or if there exist distinct
elements $x_{0},x_{1},...,x_{n}$ such that
\begin{center}
(1) $(x_{0},x_{1},...,x_{n}) \in \Delta$
\end{center}
where $x_{0}=x$ and $x_{1},...,x_{n} \in A$

and {\it directly dependent} on $(x)$ or $(x_{0},x_{1},...,x_{n})$,
respectively.
\end{defi}

\begin{defi}
A set $A$ is called {\it dependent} if (1) is satisfied for some
distinct elements $x_{0},x_{1},...,x_{n} \in A$, and otherwise $A$
is {\it independent}.
\end{defi}

\begin{defi}
If a set $A$ is  {\it independent} and for any $x \in {\bf S}$, $x
\sim \Sigma A$, i.e. $x$ is dependent on $A$, then $A$ is called a
{\it basis of} ${\bf S}$.
\end{defi}

\begin{defi}
{\bf TRANSITIVITY AXIOM:}
\begin{center}
If $x \sim \Sigma A$ and for all $a \in A$, $a \sim \Sigma B$, then
$x \sim \Sigma B$.
\end{center}
\end{defi}

Note, that the following well known properties of independent sets
are satisfied:

(2) Any subset of an independent set $A$ is independent,

(3) A basis is a maximal independent set of ${\bf S}$ and vice
versa.

(4) The family $(\mathbb{X}, \subseteq )$ of all independent subsets
of ${\bf S}$ is partially ordered by the set-theoretical inclusion.
Shortly we say that $\mathbb{X}$ is an ordered set (a po-set).

(5) Any superset of a dependent set of ${\bf S}$is dependent.

{\it Proof} (2) Let $A$ be an independent set and $B \subseteq A$.
Then $B$ is independent, as if not, then in $B$ there is a finite
sequence $b_{0},...,b_{n}$ of elements such that $(b_{0},...,b_{n})
\in \Delta$. But $(b_{0},...,b_{n}) \in A$, therefore $A$ is
dependent, a contradiction.

(3) Let $A$ be a basis of ${\bf S}$, i.e. $A$ is independent and for
each $x \in {\bf S}$, $x \sim \Sigma A$. Assume that $A$ is not a
maximal independent set of $S$, let $A \subset B$, with $A \neq B$,
where $B$ is independent and let $b \in B - A$. We have $b \sim
\Sigma A$, as $A$ is a basisis, i.e. there is a sequence of
elements: $(x_{0}=b, x_{1},...,x_{n})$ with $x_{1},...,x_{n} \in A$,
and such that $x_{0}, x_{1},...,x_{n} \in \Delta$. We get that $B$
is dependent, a contradiction. We conclude that a basis is a maximal
independent set in ${\bf S}$.

Now, let $A$ be a maximal independent subset of ${\bf S}$. We show
that $A$ is a basis. We need to show that for every $x \in {\bf S}$,
$x \sim \Sigma A$. \\If $x \in A$, then $ x \sim \Sigma A$ by
Definition. If $x$ is not included in $A$, then the set $B = A + \{
x \}$ is dependent, i.e. there is a sequence
$(x_{0},x_{1},...,x_{n})$ of elements of $B$, such that
$(x_{0},x_{1},...,x_{n}) \in \Delta$. But A is independent,
therefore one of $x_{i}$, say $x_{0}$ is in the set $B - A$, i.e.
$x_{0} = x$ and all $x_{1},...,x_{n} \in A$ (as all $x_{i}$ are
different, for $i=0,1,...,n$). We obtain $x \sim \Sigma A$.
Therefore $A$ is a basis of ${\bf S}$. $\Box$

\section{Examples}
\begin{example}
Consider the two-dimensional vector space ${\bf R}^2$.
\\$\Delta$ is defined as follows. A subset $X$ of $R^2$ containing
two parallel vectors is always dependent. A finite subset $X$ of
$R^2$ with more that 2 elements is always dependent. The one element
subset $\{ 0 \}$ is defined as dependent. Transitivity axiom is
satisfied in such dependance space ${\bf S}$. Independent subsets
contain at lest two non parallel vectors. Bases are bases
 of ${\bf R}^2$ in the classical sense.
\end{example}
\begin{example}
Let ${\bf K}$ denotes the set of some (at least two) colours, ${\bf
S} = C(K)$ is the set of all sequences of $K$. $\Delta =
C_{fin}({\bf K})$ be the set containing all finite at lest
two-element sequences of elements of ${\bf K}$ with at least one
repetition of colours. Then the Transitivity Axiom is satisfied for
such defined dependence space ${\bf S}$.
\end{example}

\begin{example}
Let $C({\bf N})$ denotes the set of all sequences of natural numbers
and $\Delta = C_{fin}({\bf N})$ be the set containing all finite at
lest two-element sequences of elements of ${\bf N}$ with at least
one repetition. Then the Transitivity Axiom is satisfied for such
defined dependence space ${\bf S}$.
\end{example}

\begin{example}
A graph is an abstract representation of a set of objects where some
pairs of the objects are connected by links. The interconnected
objects are represented by mathematical abstractions called
vertices, and the links that connect some pairs of vertices are
called edges. Typically, a graph is depicted in diagrammatic form as
a set of dots for the vertices, joined by lines or curves for the
edges. For a graph, call a finite subset $A$ of its vertices
directly dependent if it has at least two elements which are
connected. Then for a vertex $x$, put $x \sim \Sigma A$ if there
exists a link between $x$ and a vertex $a$ in $A$. Then the
transitivity axiom for such a dependence space is satisfied.
\end{example}

\section{Po-set of independent sets}

Following K. Kuratowski and A. Mostowski \cite{4} p. 241, a po-set
$(\mathbb{X}, \subseteq)$ is called {\it closed} if for every chain
of sets $\mathbb{A}  \subseteq P(\mathbb{X})$ there exists $\cup
\mathbb{A}$ in $
\mathbb{X}$, i.e. $\mathbb{A}$ has the supremum in
$(\mathbb{X}, \subseteq)$.

\begin{thm}
The po-set $(\mathbb{X}, \subseteq)$ of all independent subsets of
${\bf S}$ is closed.
\end{thm}

{\it Proof} Let $\mathbb{A}$ be a chain of independent subsets of
${\bf S}$, i.e. $\mathbb{A} \subseteq P(\mathbb{X})$, and for all
$A,B \in \mathbb{A}$ $(A \subseteq B)$ or $(B \subseteq A)$. We show
that the set $\cup\mathbb{A}$ is independent. Otherwise there exist
elements $(x_{0},x_{1},...,x_{n}) \in \Delta$ such that $x_{i} \in
\cup\mathbb{A}$, for $i=0,...,n$. therefore there exists a set $A
\in \mathbb{A}$ such that $x_{i} \in A$ for all $i=0,...,n$.
\\We conclude that $A$ is dependent, $A \in \mathbb{A}$, a
contradiction. $\Box$

\section{Steinitz' exchange theorem}

A transfinite version of the Steinitz Exchange Theorem, provides
that any independent subset injects into any generating subset. The
following is a generalization of Steinitz' Theorem proved originally
in 1913 and then in \cite{1}--\cite{2}:
\begin{thm}
If $A$ is a basis and $B$ is an independent subset (of a dependence
space ${\bf S}$). Then assuming Kuratowski-Zorn Maximum Principle,
there is a definite subset $A'$ of $A$ such that the set $B + (A -
A')$ is also a basis of ${\bf S}$.
\end{thm}

{\it Proof} If $B$ is a basis then $B$ is a maximal independent
subset of ${\bf S}$ and $A' = A$ is clear.

Assume that $A$ is a basis and $B$ is an independent subset (of the
dependence space ${\bf S}$). Consider $\mathbb{X}$ to be the family
of all independent subsets of ${\bf S}$ containing $B$ and contained
in $A+B$. then $(\mathbb{X}, \subseteq)$ is well ordered and closed.
Therefore assuming Kuratowski-Zorn Maximal Principle \cite{3} there
exists a maximal element of $\mathbb{X}$. \\We show that this
maximal element $X \in \mathbb{X}$ is a basis of ${\bf S}$.

As $X \in \mathbb{X}$ then $B \subseteq X \subseteq A+B$ by the
construction. Therefore $X = B + (A - A')$ for some $A' \subseteq
A$. We show first that for all $a \in A$, $a \sim \Sigma X$. If $a
\in X$ then $a \sim \Sigma X$ by the definition. If not, then put $Y
= X + \{ a \}$. Then $X \neq Y$, $X \subseteq Y$, $B \subseteq Y
\subseteq A + B$ and $Y$ is dependent in ${\bf S}$.

By the definition there exist elements: $(x_{0},x_{1},...,x_{n}) \in
\Delta$ with $x_{1},...,x_{n} \in X + \{ a \}$, as $X$ is an
independent set. Moreover, one of $x_{i}$ is $a$, say $x_{0} = a$.
We get $a \sim \Sigma X$ as $x_{0},x_{1},...,x_{n}$ are different.

Now we show that $X$ is a basis of ${\bf S}$. Let $x \in S$, then $x
\sim \Sigma A$ as $A$ is a basis of ${\bf S}$. Moreover for all $a
\in A$, $a \sim \Sigma X$, thus $x \sim \Sigma X$ by the
Transitivity Axiom. $\Box$

\end{document}